\documentclass[12pt]{amsart}
 \usepackage{amsmath}
 \usepackage{amssymb}
 \usepackage{amscd}
 \usepackage[mathscr]{eucal}
\usepackage{epsfig}
\usepackage{graphics}

 \newtheorem{thm}{Theorem}[section]
 \newtheorem{lem}[thm]{Lemma}
 
 \newtheorem{example}[thm]{Example}
 \newtheorem{cor}[thm]{Corollary}
 \newtheorem{prop}[thm]{Proposition}
 \newtheorem{defn}[thm]{Definition}

 \newtheorem{prob}[thm]{Problem}
 \newcommand{\Hom}{\operatorname{Hom}}
\newcommand{\lk}{\operatorname{lk}}
\newcommand{\inv}{\operatorname{inv}}
\newcommand{\des}{\operatorname{des}}

\newcommand{\A}{\mathscr{A}}
\newcommand{\M}{\mathscr{M}}
\renewcommand{\SS}{\mathscr{S}}

\newcommand{\Z}{\mathbb{Z}}
\newcommand{\N}{\mathbb{N}}
\newcommand{\R}{\mathbb{R}}
\newcommand{\Th}{\tilde{h}}

\begin{document}
 
 \title[$g$-elements]{$g$-elements, finite buildings and higher Cohen-Macaulay connectivity}
\thanks{Partially supported by NSF grant DMS-0245623.}
\author{Ed Swartz}
\email{ebs22@cornell.edu}
\address{Dept. of Mathematics, Cornell University, Ithaca, NY, 14853.}

\begin{abstract}
 Chari proved that if $\Delta$ is a $(d-1)$-dimensional simplicial complex with  a convex ear decomposition, then $h_0 \le \dots \le h_{\lfloor d/2 \rfloor}$ \cite{Ch}. Nyman and Swartz raised the problem of whether or not the corresponding $g$-vector is an $M$-vector  \cite{NS}.  This is proved to be true by showing that the set of pairs $(\omega,\Theta),$ where $\Theta$ is a l.s.o.p. for $k[\Delta],$ the face ring of $\Delta,$ and $\omega$ is a $g$-element for $k[\Delta]/\Theta,$ is nonempty whenever the characteristic of $k$ is zero.

Finite buildings have a convex ear decomposition.  These decompositions point to inequalities on the flag $h$-vector  of such spaces similar in spirit to those examined in \cite{NS} for order complexes of geometric lattices.   This also leads to connections between higher Cohen-Macaulay connectivity and conditions which insure that $h_0 < \dots < h_i$ for a predetermined $i.$

 \end{abstract}

  \maketitle
  
One of the most basic combinatorial invariants of a (finite) simplicial complex is its $f$-vector, or equivalently, its $h$-vector.  In order to analyze $h$-vectors of matroid independence complexes Chari introduced the notion of a convex ear decomposition \cite{Ch}.  He showed that $(d-1)$-dimensional  complexes which have such a decomposition satisfy $h_i \le h_{d-i}$ and $h_i \le h_{i+1}$ for all $i \le \lfloor d/2 \rfloor.$  In addition, he proved that independence complexes of matroids have a PS-ear decomposition, a special type of convex ear decomposition.  Spaces with a PS-ear decomposition satisfy the additional condition that their $g$-vector, $(g_0,g_1, \dots, g_{\lceil d/2 \rceil})$, where $g_i = h_i - h_{i-1},$ is an M-vector \cite{Sw3}.  Our main result, Theorem \ref{main}, is that this holds for all spaces with a convex ear decomposition.

In section \ref{buildings} we introduce a convex ear decomposition for finite buildings.  In addition to the enumerative conclusions above, this will allow an analysis of the flag $h$-vectors of such complexes.   We end with an examination of a connection between higher Cohen-Macaulay connectivity and increasing $h$-vectors.   

 Throughout, $\Delta$ is a finite $(d-1)$-dimensional abstract simplicial complex with vertex set $V, |V|=n.$ A maximal face (under inclusion) of $\Delta$ is a {\it facet.} The {\it $f$-vector} of $\Delta$ is $(f_0,\dots, f_d),$ where $f_i$ is the number of faces of $\Delta$ of cardinality $i.$ (Note: Our $f_i$ is frequently denoted by $f_{i-1}.$) The {\it $h$-vector} of $\Delta$ is $(h_0,\dots,h_d)$ where
 
 \begin{equation} \label{h by f}
  h_i(\Delta) = \sum^i_{j=0} (-1)^{i-j} \binom{d-j}{d-i} f_j(\Delta).
\end{equation}  
 
 \noindent Equivalently,
 
 \begin{equation} \label{f by h}
  f_j(\Delta) = \sum^j_{i=0} \binom{d-i}{d-j} h_i(\Delta).
\end{equation} 

 We use $\Delta - v$ for the complex consisting of $\Delta$ with all of the faces containing the vertex $v$ removed.  Similarly, if $A \subseteq V,$ then $\Delta-A$ is $\Delta$ with all of the faces which contain any vertex in $A$ deleted.  

The {\it order complex} of a poset $P$ is the simplicial complex whose faces are the chains in $P.$  However, if $P$ contains a maximal element $\hat{1}$ or a minimal element $\hat{0},$ then we will always assume that the order complex refers to the poset $P - \{\hat{1}, \hat{0}\}.$

\begin{section}{Convex ear decompositions} \label{convex ears}

A {\it convex ear decomposition} of $\Delta$ is an
ordered sequence $\Delta_1,\dots,\Delta_m$ of
pure $(d-1)$-dimensional subcomplexes of $\Delta$ such that 

\begin{enumerate}
 \item
  $\Delta_1$ is the boundary complex of a simplicial
$d$-polytope. For each $j=2,\dots,m,
\Delta_j$ is a $(d-1)$-ball which is a proper subcomplex
of the boundary of a simplicial $d$-polytope.

  \item
    For $j \ge 2, \Delta_j \cap (\bigcup ^{j-1}_{i=1}
\Delta_i) = \partial \Delta_j.$
 \item
$\bigcup^{m}_{j=1} \Delta_j = \Delta$.

\end{enumerate}

\noindent  The {\it initial} subcomplex is $\Delta_1.$ Each $\Delta_j,$ for $j \ge 2,$ is an {\it ear} of the decomposition.  

Convex ear decompositions were originally introduced by Chari.  His original example of a convex ear decomposition was the independence complex of a matroid.  In fact, he proved that the independence complex of a matroid has a special type of convex ear decomposition, a {\it PS-ear decomposition.}  In a PS-ear decomposition the initial subcomplex is a join of boundaries of simplices, and each ear is a join of a simplex and boundaries of simplices \cite{Ch}.  

Using an idea of Bj\"{o}rner \cite{Bj2}, Nyman and Swartz showed that order complexes of geometric lattices have a convex ear decomposition \cite{NS}. In this case, the  initial complex is the first barycentric subdivision of the boundary of a simplex and each ear is a shellable ball which is a  subcomplex of such a space.  In addition to the enumerative conclusions of Theorem \ref{main} below, this approach led  to several inequalities for the flag $h$-vector of such complexes.   

\begin{defn}
  A {\bf balanced} complex is a $(d-1)$-dimensional simplicial complex $\Delta$ and a map $\phi: V \to S, \ |S|=d,$ such that $\phi(v) \neq \phi(w)$ for any pair of distinct vertices $v$ and $w$ which are contained in a face of $\Delta.$ 
\end{defn}

 Equivalently, the one-skeleton of $\Delta$ is properly $d$-colorable.  Our balanced complexes were called completely balanced in \cite{St5}.  A common example of a  balanced complex is the order complex of a ranked poset, with $\phi(x)$ the rank of $x.$ For balanced complexes there is a natural refinement of the $f$- and $h$-vectors.

\begin{defn}  Let $(\Delta, \phi)$ be a balanced complex.  Let $A \subseteq S.$  
$$\begin{array}{lcl}
\Delta_A  & = & \{\rho \in \Delta: \forall \ v \in \rho, \phi(v) \in A\}.\\
f_A &  =  & f_{|A|}(\Delta_A). \\
h_A &  = & \displaystyle\sum_{B \subseteq A} (-1)^{|A-B|} f_B.
\end{array}$$
\end{defn}

\begin{prop}\cite[pg. 95--96]{St}
  Let $(\Delta,\phi)$ be a pure balanced complex.  Then
  $$f_i = \displaystyle\sum_{|A|=i} f_A, $$
  $$h_i = \displaystyle\sum_{|A|=i} h_A.$$
\end{prop}

The inequalities for the flag $h$-vector of the order complex of a rank $d+1$ geometric lattice in \cite{NS} can be described in terms of the weak  order (also known as the weak Bruhat order) on the symmetric group $S_{d+1}.$ Let $\pi$ be permutation in $S_{d+1}.$  The {\it inversion} set of $\pi$ is the set of all pairs $\inv (\pi) = \{(i,j): 1 \le i < j \le d+1, \pi^{-1}(i) > \pi^{-1}(j)\}.$ 

\begin{defn}  The {\bf weak } order on $S_{d+1}$ is defined by
$$ \pi \le  \pi^\prime \leftrightarrow \inv(\pi) \subseteq \inv (\pi^\prime).$$
\end{defn}
  
  The {\it descent set} of $\pi \in S_{d+1}$ is $\des (\pi) = \{i: \pi(i) > \pi(i+1)\}.$  Obviously, $\des(\pi) \subseteq [d].$ Given $A \subseteq [d],$ the {\it descent class} of $A$ is $D(A) = \{\pi \in S_{d+1}: \des(\pi) = A\}.$  Finally, for $A$ and $B$ subsets of $[d]$ we say $A$ {\it dominates} $B$ if there exists an injection $\psi:D(B) \hookrightarrow D(A)$ such that $\pi \le  \psi(\pi)$ for all $\pi \in D(B).$
  
  \begin{thm} \cite{NS} \label{geom lattice}
    Let $\Delta$ be the order complex of a rank $d+1$ geometric lattice.  If $A$ and $B$ are subsets of $[d]$ such that $A$ dominates $B,$ then
    $$h_B \le h_A.$$
  \end{thm}

Several  order complexes of  posets have convex ear decompositions.  These include  rank-selected  subposets of geometric lattices, supersolvable  lattices with nonzero M\"{o}bius function on every interval and their rank-selected subposets,  and $d$-divisible partition lattices ($d \ge 3$) \cite{Schw}.  The flag $h$-vectors of rank $(d+1)$ supersolvable lattices with nonzero M\"obius function on every interval also satisfy the conclusion of Theorem \ref{geom lattice}.

There is a general construction which includes all of the above examples and the buildings in the next section.  Let $\Sigma$ be a contractible, shellable $d$-polytopal complex. For polytopal shellings see, for instance, \cite[chapter 8]{Zi}.  Let $\Sigma_1, \dots, \Sigma_m$ be a shelling order of the facets and  assume that all of the facets are  simplicial $d$-polytopes.  Removing all the open $d$-cells leaves a $(d-1)$-dimensional simplicial complex with a convex ear decomposition.  The initial subcomplex is $\partial (\Sigma_1).$  For $2 \le j \le m,$ the ear $\Delta_j$ is the closure of $\partial (\Sigma_j) - \cup^{j-1}_{i=1} \partial (\Sigma_i).$

\end{section} 
  
 \begin{section}{finite buildings} \label{buildings} 
  
  Finite buildings   have a convex ear decomposition.  This is an immediate consequence of \cite[Lemma 3.5]{He}. We will use this decomposition when examining the complementary $h$-vector in Section \ref{kcm}.  However, for reasons we will make clear below (see Theorem \ref{building flag h}), we prefer another proof here.  

There are several standard references on buildings.  We mention \cite{Bj3}, \cite{Br}, and \cite{Ti}, as all of the facts we use can be found in those references.   Let $(W,S)$ be a finite Coxeter system with associated Coxeter complex $\Sigma(W,S).$ Specifically, $W$ is a finite group generated by reflections of (linear) hyperplanes  in $U,$ a $d$-dimensional real vector space, and $S$ is a generating set of reflections defined below.  The collection of hyperplanes is assumed to be essential, that is their intersection is the origin,  and contains all of the  hyperplanes of the reflections in $W.$   The intersection of the unit sphere of $U$ with the hyperplane arrangement results in a (spherical) simplicial complex $\Sigma(W,S)$ which is a triangulation of the $(d-1)$-sphere.   

The group $W$ acts transitively and freely on the facets (also called chambers) of $\Sigma(W,S).$  Let $\sigma$ be a fixed facet of $\Sigma(W,S).$ The simply transitive action of $W$ on the facets allows us to  identify $w \in W$ with the facet $w \cdot \sigma.$  The linear span of each $(d-2)$-face of the boundary of $\sigma$ is one of the  hyperplanes in the arrangement. The corresponding set of reflections is  $S = \{s_1,\dots,s_{d-1}\}$ and generates $W$. 

Given $w \in W$, the minimal $\ell$ such that $w = s_{i_1} \cdots s_{i_\ell}$ is  $\ell(w),$ the {\it length} of $w.$  
The {\it weak } order on $W$ (also known as the weak Bruhat order) is defined by $w <  w^\prime$ if there exists $s_1,\dots,s_j \in S$ such that $w \cdot s_1 \cdots s_j = w^\prime$ and $\ell(w^\prime) = \ell(w) + j.$ An equivalent formulation given by the geometry of $\Sigma(W,S)$ is as follows.  A {\it path} (also called a gallery) in $\Sigma(W,S)$ is a sequence, $(\sigma_0,\dots,\sigma_t),$ of facets such that for each $i$ the intersection of $\sigma_i$ and $\sigma_{i+1}$ is a $(d-2)$-face (usually called a wall).  The length of the path is $t.$  A {\it geodesic} is a path of minimal length among all paths beginning and ending with the same facets.  The weak  order on $(W,S)$  is equivalent to $w \le  w^\prime$ if and only if there is a geodesic from $\sigma$ to $w^\prime \cdot \sigma$ which contains $w \cdot \sigma.$ 

\begin{example}
  Let $W$ be the group generated by the collection of hyperplanes $H_{ij}=\{(x_1,\dots,x_{d+1}) \in \R^{d+1}: x_i = x_j\}, i< j.$  This is not an essential arrangement as it contains the line $x_1= \dots = x_{d+1}.$  Intersecting with the $d$-dimensional subspace $U=\{(x_1,\dots,x_{d+1}) \in \R^{d+1}: x_1 + \dots + x_{d+1} = 0\}$  does give an essential arrangement.  The group generated by the corresponding reflections in $U$  is the symmetric group $S_{d+1}.$  If we choose $\sigma$ to be the facet which contains   those $(x_1,\dots,x_{d+1})$ such that $x_1 < \dots < x_{d+1},$ then $S$ corresponds to the transpositions $(i,i+1)$ and the weak  order on $(W,S)$ is the same as defined in Section \ref{convex ears}.
  \end{example}

The {\it descent} set of $w \in W$ is $\des (w) = \{ s \in S: \ell(w \cdot s) < \ell(w)\}.$  The descent set of $w$ can also be defined via $\Sigma(W,S).$ 
Let $V_{(W,S)}$ be the vertices of the Coxeter complex $\Sigma(W,S).$ For each $(d-2)$-face $\tau$ of $\partial \sigma,$ set $\psi(\tau)$ to be the corresponding hyperplane in $S.$  Define $\psi: V_{(W,S)} \to S$ by first defining $\psi$ on the vertices of $\sigma$ to be $\psi$ of the opposite face.  Then extend $\psi$ to a labeling of all of $V_{\Sigma(W,S)}$ in the only way possible that insures that $\Sigma(W,S)$ and $\psi$ form a balanced complex.  This also labels the $(d-2)$-faces of $\Sigma(W,S).$ Simply assign $\psi(\tau)$ to be $\psi(v),$ where $v$ is any vertex such that $\tau \cup \{v\}$ is a facet.   With this definition of $\psi$ on the $(d-2)$-faces of $\Sigma(W,S),$ the descent set of $w$ is the set of all $s$ such that there exists a geodesic $(\sigma, \dots, \sigma^\prime, w \cdot \sigma)$ with $\psi(\sigma^\prime \cap (w \cdot \sigma)) = s.$

\begin{defn}
  Let $\Sigma(W,S)$ be a finite Coxeter complex.  A {\bf finite building} of type $(W,S)$ is a (finite) simplicial complex $\Delta$ which is  the union of subcomplexes $\Sigma,$ called {\bf apartments}, such that
  
  \begin{itemize}
    \item
      Each apartment $\Sigma$ is isomorphic to $\Sigma(W,S).$
    \item
      For any two faces $\rho_1$ and $\rho_2$ in $\Delta,$ there is an apartment $\Sigma$ containing both of them.
    \item
      If $\Sigma$ and $\Sigma^\prime$ are two apartments containing $\rho_1$ and $\rho_2$, then there is an isomorphism $\Sigma \to \Sigma^\prime$ fixing $\rho_1$ and $\rho_2$ pointwise.
    \end{itemize}
  \end{defn}

 For the rest of this section  $\Delta$ is a finite building of type $(W,S).$ Let $\tau$ be a facet of $\Delta$ and let $\rho$ be any face of $\Delta.$  A geodesic from $\tau$ to $\rho$ is a geodesic $(\tau=\tau_0,\dots,\tau_t)$ such that $\rho \not\subseteq \tau_i$ for any $i < t$ and $\rho \subseteq \tau_t.$ There exists a unique facet $p_\rho(\tau),$ the {\it projection} of $\tau$ on $\rho,$ such that every geodesic from $\tau$ to $\rho$ ends with $p_\rho(\tau)$  \cite[3.18 - 3.19]{Ti}.
  
  Let $\sigma$ be any fixed base facet.  A facet $\sigma_i$ is {\it opposite} $\sigma$ if it is maximally distant from $\sigma.$  Let $\sigma_1, \dots, \sigma_m$ be the facets opposite  $\sigma.$ For each $\sigma_j,$  there is a  unique apartment, $\Sigma_j,$ which contains $\sigma$ and $\sigma_j.$  It is the union of all geodesics from $\sigma$ to $\sigma_j.$  Finally, set $\Delta_1 = \Sigma_1$ and for $j \ge 2$, define $\Delta_j$ to be the union of the facets of $\Sigma_j$ not contained in any $\Sigma_i, i < j.$  Since $\Delta = \cup \Sigma_j,$ it is immediate that $\Delta = \cup  \Delta_j.$

\begin{thm} \label{building ears}
  Let $\Delta, \Delta_1, \dots, \Delta_m$ be as above.  Then $\Delta_1, \dots, \Delta_m$ is a convex ear decomposition of $\Delta.$
\end{thm}

\begin{proof}
We begin by proving  that for $j \ge 2$, $\Delta_j$ is a shellable subcomplex of $\Sigma_j,$ and hence a ball.   Since $\sigma_j$ is not in any other $\Sigma_i, \ \sigma_j \in \Delta_j.$  What other facets of $\Sigma_j$ are contained in $\Delta_j$?  A facet $\tau \in \Sigma_j$ is in $\Delta_j$ if and only if for all $i < j, \ \tau$ it is not contained in any geodesic from $\sigma$ to $\sigma_i.$

 Let $\tau_1, \dots, \tau_t$ be an ordering of the facets in $\Delta_j$ which is a linear extension of the order dual of the  weak  order restricted to $\Delta_j.$   Specifically, if $\tau <  \tau^\prime,\tau = \tau_k,$ and  $\tau^\prime = \tau_l,$ then $l< k.$  From the above discussion we know that if $\tau <  \tau^\prime$ and $\tau \in \Delta_j,$ then $\tau^\prime \in \Delta_j.$  Thus, $\tau_1, \dots, \tau_t$ is an initial segment of a linear extension of the order dual of the weak  order on all of $\Sigma_j.$  By \cite[Thm A.1]{Bj3}, $\tau_1, \dots, \tau_t$  is an initial segment of a shelling of $\Sigma_j$ and hence $\Delta_j$ is a ball.
 
 In order to see that $\Delta_j \cap \cup^{j-1}_{i=1} \Delta_i = \partial \Delta_j,$ we first note that a face  $\rho$ is in $\partial \Delta_j$ if and only if $\rho$ is contained in facets  $\tau_1$ and $\tau_2$ with $\tau_1 \in \Delta_j$ and $\tau_2 \in \Sigma_j - \Delta_j.$ Now suppose $\rho$ is  a face in $\Delta_j$ and $\Delta_i$ with $i<j.$  Let $\tau=p_\rho(\sigma_i),$ the projection of $\tau$ on $\rho.$  Since $\rho \in \Sigma_i \cap \Sigma_j,  \ \tau$ must also be in $\Sigma_i \cap \Sigma_j.$  Hence, $\rho \subseteq \tau$ with $\tau$  a facet of $\Sigma_j$ not in $\Delta_j.$  Thus $\rho \subseteq \partial \Delta_j.$  The other inclusion is obvious.  
 
   \end{proof}

The above construction suggests inequalities for the flag $h$-vector of a finite building similar in spirit to Theorem \ref{geom lattice}.  For $A \subseteq S$ the {\it descent class} of $A$ is $D(A) = \{ w \in W: \des w = A\}.$  If $A,B \subseteq S,$ then $A$ {\it dominates} $B$ if there exists an injection $\psi: D(B) \hookrightarrow D(A)$ such that $w <  \psi(w)$ for all $w \in D(B).$

\begin{thm}  \label{building flag h}
  Let $(W,S)$ be a finite Coxeter system.  Let  $A$ and $B$ be subsets of $S$ and assume $A$ dominates $B.$  If $\Delta$ is a finite building of type  $(W,S),$ then 
  
  $$h_B \le h_A.$$
  
  \end{thm}    
  
  \begin{proof}
    Let $\psi:D(B) \hookrightarrow D(A)$ be an injection such that $w <  \psi(w)$ for all $w \in D(B).$  Let $\sigma$ be a fixed base facet. For any $B \subseteq S, h_B$ is the number of facets $\tau$ in $\Delta$ whose descent set is $B$, where the descent set is computed in any apartment containing $\sigma$ and $\tau$ \cite{Bj3}. By the proof of Theorem \ref{building ears}, each $\tau$ with descent set $B$ is in exactly one $\Delta_j,$ and if $\tau$ corresponds to $w$ in $\Sigma_j,$ then the facet $\tau^\prime$ which corresponds to $\psi(w)$ in $\Sigma_j$ is also in $\Delta_j$ and has descent set $A.$  Hence, there is an injection from facets with descent set $B$ to facets with descent set $A.$
    
    \end{proof} 
  
  For buildings associated to $GF(q),$ the above theorem follows easily from the formula \cite{Bj3}
  
 \begin{equation} \label{GF(q) flagh}
 h_A = \displaystyle\sum_{w \in D(A)} q^{\ell(w)}.
 \end{equation}

    One possible approach to looking for a combinatorial proof of some of the inequalities implied by Corollary \ref{main enumerate} for finite buildings would be to consider the following problems.
\begin{prob}
Let  $(W,S)$ be a finite Coxeter system and let $d-1$ be the dimension of the Coxeter complex associated to $(W,S).$
  \begin{enumerate}
    \item
      For which subsets $A,B$ of $S$ does $A$ dominate $B$?  
    \item Suppose $i \le d/2.$ Does there exist an injection $\iota$ from $i$-subsets of $S$ to $(i+1)$-subsets of $S$ such that $\iota(B)$ dominates $B$ for all $i$-subsets $B$?  
    \item Is there an injection $\alpha: \displaystyle\cup_{|B|=i} D(B) \hookrightarrow \displaystyle\cup_{|A|=i+1} D(A)$  such that $w <  \alpha(w)$?
    \item Is there a bijection $\beta:\displaystyle\cup_{|B|=i} D(B) \hookrightarrow \displaystyle\cup_{|A|=d-i} D(A)$  such that $w <  \beta(w)$?
    \end{enumerate}
  \end{prob}

\noindent Some of these problems were explored for the symmetric group in \cite{NS}.

\end{section}

\begin{section}{Face rings} \label{face rings}

One of the most powerful tools for studying enumerative properties of simplicial complexes, especially Cohen-Macaulay complexes, is the {\it face ring}, also known as the {\it Stanley-Reisner ring}.   Let $k$ be any field and set $R = k[x_1,\dots,x_n].$ For any homogeneous ideal $I$ of $R$ we use $(R/I)_i$ to represent the degree $i$ component of $R/I.$

\begin{defn}
  The {\bf face ring} of $\Delta$ is 
$$k[\Delta] = R/I_\Delta,$$
\noindent where $I_\Delta = <\{x_{i_1} \cdots x_{i_j}: \{v_{i_1},\dots, v_{i_j}\} \notin \Delta\}>.$
\end{defn}

Let $\Theta=\{\theta_1,\dots,\theta_d\}$ be a set of one-forms in $R.$ We will also use $\Theta$ to represent the ideal $<\Theta>$ in $R$ or $k[\Delta],$ relying on context to make it clear what is intended.  Write each $\theta_i = \theta_{i1} x_1 + \dots + \theta_{in} x_n$ and let $ T$ be the matrix $(\theta_{ij}).$  
To each facet $\sigma$ in $\Delta$ let $T|_\sigma$ be the submatrix of $T, (\theta_{ij})_{v_j \in \sigma},$ consisting of the columns of $T$ corresponding  to the vertices of $\sigma.$   If $T|_\sigma$ has rank $|\sigma|$ for every facet $\sigma$, then $\Theta$ is a {\it linear set of parameters} 
(l.s.o.p.) 
for $k[\Delta].$  If $k$ is infinite, then it is always possible to choose $\Theta$ 
such that every set 
of $d$ columns of $T$ is independent. 

Two of the most useful facts concerning $k[\Delta]$ are the following.

\begin{thm} \cite{Re}  Let $k[\Delta]$ be the face ring of $\Delta.$  Then
  $k[\Delta]$ is a Cohen-Macaulay ring if and only if for all faces $\sigma,$
$$\tilde{H}_j(\lk \sigma; k) = 0 \mbox{ for all } j < d-|\sigma|-1.$$
\end{thm} 

 We call $\Delta$ a {\it Cohen-Macaulay} complex if $k[\Delta]$ is a Cohen-Macaulay ring.  It can be shown that the property of being Cohen-Macaulay is a purely topological property \cite{Mu}.  

A Mayer-Vietoris argument applied inductively to the number of ears shows that spaces with a convex ear decomposition are Cohen-Macaulay.  As we shall see in Section \ref{kcm},  they are doubly Cohen-Macaulay.

\begin{thm} \cite{St1}  \label{CM}
If $\Delta$ is a Cohen-Macaulay complex, then for any l.s.o.p. $\Theta,$
$$h_i(\Delta) = \dim_k (k[\Delta]/\Theta)_i.$$
\end{thm}

An immediate consequence of this result is that $h_i \ge 0$ for any Cohen-Macaulay complex.  In addition,  the $h$-vector of a CM complex is an M-vector.  A sequence $(h_0,h_1,\dots,h_d)$ is an {\it M-vector} if it is the Hilbert function of a homogeneous quotient of a polynomial ring. A purely arithmetic criterion is the following description due  to Macaulay.

Given  positive integers $h$ and $i$
there is a  unique way of writing
$$h=\binom{a_i}{i} + \binom{a_{i-1}}{i-1} + \dots + \binom{a_j}{j}$$
so that $a_i > a_{i-1} > \dots > a_j \ge j \ge 1.$  Define
$$h^{<i>} = \binom{a_i+1}{i+1} + \binom{a_{i-1}+1}{i} + \dots + 
\binom{a_j+1}{j+1}$$

\begin{thm} \cite[Theorem 2.2]{St} \label{CM h-vectors}
  A sequence of nonnegative integers $(h_0, \dots, h_d)$ is an M-vector if and only if $h_0=1$ and 
  $h_{i+1} \le h^{<i>}_i$ for all $1 \le i \le d-1.$
\end{thm}

  We denote the  {\it canonical module} of $k[\Delta]$ by $\Omega(k[\Delta]).$  The canonical module is an $R$-module and can be defined using homological methods.  In fact, $I_\Delta$ is contained in the annihilator of  $\Omega(k[\Delta]),$ so the canonical module   is also a  $k[\Delta]$-module.  The only properties of  $\Omega(k[\Delta])$ that we will use are in the following theorem.  As usual, $k[\Theta]$ is the ring $k[\theta_1,\dots,\theta_d].$

\begin{thm}\cite[sections I.12, II.7]{St} \label{canonical}
 Let $\Delta$ be a $(d-1)$-dimensional Cohen-Macaulay complex and let $\Theta = \{\theta_1,\dots,\theta_d\}$ be a l.s.o.p. for $k[\Delta].$
  \begin{enumerate}
    \item
       $\Omega(k[\Delta]) \cong \Hom_{k[\Theta]}(k[\Delta],k[\Theta])$, where the $R$-module structure is given by $(f\phi)(p) = \phi(f \cdot p), f \in R, \phi \in \Omega(k[\Delta])$ and $p \in k[\Delta].$
    \item \label{grading}
      There is a grading of  $\Omega(k[\Delta])$ so that as an $\N$-graded $R$-module $$\dim_k  (\Omega(k[\Delta])/\Theta \Omega(k[\Delta]))_i = \dim_k (k[\Delta]/\Theta)_{d-i}.$$

\end{enumerate}
\end{thm}

\noindent   The grading of $\Omega(k[\Delta])$ is such that if $f \in \Omega(k[\Delta])$ maps degree $i$ elements of $k[\Delta]$ to degree $i-d$ elements of $k[\Theta],$ then $f$ has degree $0.$  Note that this is a shift of the usual grading of $\Omega(k[\Delta])/ \Theta \Omega(k[\Delta])$ whose first nonzero module is in degree $-d$ instead of $0.$

Let $\Delta$ be a Cohen-Macaulay complex and let $\Theta$ be a l.s.o.p. for $k[\Delta].$  A {\it $g$-element} for $k(\Delta)=k[\Delta]/\Theta$ is a one-form $\omega \in R$ such that  multiplication
$$\omega^{d-2i}: k(\Delta)_i \to k(\Delta)_{d-i}$$ 

\noindent is an injection for every $i, 0 \le i \le \lfloor d/2 \rfloor.$ When the multiplication maps are isomorphisms $\omega$ is usually called a Lefschetz element.

 Let $G(\Delta)$ be the set of all pairs $(\omega,\Theta) \subseteq k^{n(d+1)}$ such that $\Theta$ is a l.s.o.p. for $k[\Delta]$ and $\omega$ is a $g$-element for $k(\Delta).$  While the following is well known, we include it for the sake of completeness.

\begin{prop}
 If $\Delta$ is a Cohen-Macaulay complex, then $G(\Delta)$ is a Zariski open set.
\end{prop}

\begin{proof}
Given $i < d/2,$ let $G_i$ be the set of pairs $(\omega,\Theta)$ such that  $\omega^{d-2i}: k(\Delta)_i \to k(\Delta)_{d-i}$ is an injection and $\Theta$ is a l.s.o.p. for $k[\Delta].$  As $G(\Delta)$ is the intersection of all of the $G_i,$ it is sufficient to show that each $G_i$ is a Zariski open set.

  Since $\Delta$ is Cohen-Macaulay, it is a pure complex, so for each facet $\sigma \in \Delta, \ T|_\sigma$ is a square matrix.  Let $f_\sigma$ be its determinant.  Each $f_\sigma$ is polynomial in the $\theta_{ij}.$ Therefore, $\Theta$ is a l.s.o.p. for $k[\Delta]$ if and only if  the product of all the $f_\sigma$ is nonzero.   Denote this product by $f_\Delta.$  

Let $\M$ be the collection  of all subsets $U$ of  monomials in $R$ of degree less than or equal to $d$ such that the number of monomials of degree $i$  in $U$  is $h_i.$   The monomials in $U$ form a $k$-basis of $k(\Delta)$ if and only  if  for each $j$ the collection of all the monomials of the form $u_s \cdot v_t,$ with $u_s$ a monomial of degree $s$ in  $k[\Theta], v_t$  a monomial of degree $t$ in $U,$ and $s+t=j,$ form a basis of $k[\Delta]_j.$  Hence, there is a polynomial, $f_U,$ a product of $d+1$ determinants (one for each degree) in the $\theta_{ij}$-variables,  such that $f_\Delta f_U$  is nonzero if and only if $\Theta$ is a l.s.o.p. for $k[\Delta]$  and the monomials in $U$ are a $k$-basis of $k(\Delta).$  

Fix $U \in \M.$  For each $j$, let $U_j$ be the degree $j$ monomials in $U.$  Now we attempt to compute the matrix for multiplication  $\omega^{d-2i}:k(\Delta)_i \to k(\Delta)_{d-i}$ using the ``basis'' $U_i$  for $k(\Delta)_i,$ and the ``basis''  $U_{d-i}$ for $k(\Delta)_{d-i}.$  If  $\Theta$ is a l.s.o.p. for $k[\Delta]$  and $U_{d-i}$ is a basis for $k(\Delta)_{d-i},$ then we could compute the coefficients of $\omega^{d-2i} \cdot u$ for each $u \in U_i$ in the $U_{d-i}$-basis using Cramer's rule.  These coefficients are rational functions in the $\omega$ and $\Theta$ variables, with the denominator equal to the determinant which indicates whether or not $U_{d-i}$ is a basis of $k(\Delta)_{d-i}.$ Instead, use Cramer's rule without the divisor.  When $U_{d-i}$ is a basis for $k(\Delta)_{d-i}$ the matrix for the linear transformation $\omega^{d-2i}: k(\Delta)_i \to k(\Delta)_{d-i}$  we obtain will be a nonzero scalar multiple of the correct matrix.  In any event, the coefficients of the matrix are polynomials in the $\omega$ and $\Theta$ variables.

Let $A$ be a subset of $U_{d-i}$ of cardinality $h_i.$  If $h_i > h_{d-i}$, then there are no such $A.$  However, this implies that $G_i = \emptyset,$ a Zariski open set.  So we may assume that there are such sets.  Let $f_A$ be the determinant of the corresponding $h_i \times h_i$ minor of the matrix determined by our modified Cramer's rule.  If $f_A f_U f_\Delta (\omega, \Theta) \neq 0$ for some $U$ and $A,$  then $(\omega, \Theta) \in G_i.$  Conversely, suppose  $(\omega, \Theta) \in G_i.$ Since there exists some $U \in \M$ such that the monomials in $U$ form a $k$-basis of $k(\Delta),$  there must be some $h_i$-subset $A$ of some $U_{d-i}$ so that  $f_A f_U f_{\Delta} (\omega,\Theta) \neq 0.$   As there are only finitely many polynomials $f_A f_U f_\Delta, \ G_i$ is a Zariski open set.

\end{proof}

Using the above proposition, the necessary part of the $g$-theorem for simplicial polytopes can be stated in the following form.

\begin{thm}\cite{St2}, \cite{Mc2} \label{g thm}
Let $\Delta$ be the boundary complex of a simplicial polytope and let $k$ be a field of characteristic zero.  Then $G(\Delta)$ is not empty.
\end{thm}

Kalai and Stanley used Theorem \ref{g thm} to establish restrictions on the $h$-vectors of balls which were full dimensional subcomplexes of the boundary of a simplicial polytope.

\begin{thm}\cite{Kal},\cite{St3}  \label{surj balls}
  Suppose $\Delta$ is homeomophic to a $(d-1)$-ball and is a subcomplex of $\Sigma,$ where $\Sigma$ is the boundary of a simplicial $d$-polytope. Then for any $\omega \in G(\Sigma),$ multiplication $\omega^{d-2i}:k(\Delta)_i \to k(\Delta)_{d-i}$ is surjective.  
\end{thm}

\begin{thm} \label{main}
  If $\Delta$ has a convex ear decomposition and the characteristic of $k$ is zero,  then $G(\Delta)$ is not empty.
\end{thm}

\begin{proof}
  The proof is by induction on $m,$ the number of ears.  Theorem \ref{g thm} is $m=1.$ Let $\Sigma = \bigcup^{m-1}_{j=1} \Delta_j.$  Let $I$ be the kernel in the short exact sequence of $R$-modules
$$0 \to I \to k[\Delta] \to k[\Sigma] \to 0.$$

\noindent Evidently $I$ is the ideal of $k[\Delta]$ generated by the interior faces of $\Delta_m.$  As an $R$-module, $I$ is also the kernel in the short exact sequence
$$ 0 \to I \to k[\Delta_m] \to k[\partial \Delta_m] \to 0.$$ By a theorem of Hochster  \cite[Theorem II.7.3]{St}, $I$ is isomorphic to $\Omega(k[\Delta_m])$  as a $\Z$-graded module.

Dividing out by $\Theta$ gives the short exact sequence of $R$-modules

$$ 0 \to I/(I \cap \Theta) \to k(\Delta) \to k(\Sigma) \to 0.$$

\noindent  For each $i, \dim_k (I/(I \cap \Theta))_i = h_i(\Delta) - h_i(\Sigma)=h_i(\Delta_m) - g_i(\partial \Delta_m).$ This is $h_{d-i}(\Delta_m)$ \cite{St3}.  Now,  $(I/\Theta I)_i \cong (\Omega(k[\Delta_m])/ \Theta \Omega(k[\Delta_m]))_i,$ so $\dim_k (I/\Theta I)_i = h_{d-i}(\Delta_m).$  Since $\Theta I \subseteq I \cap \Theta$ we must have 

$$\Omega(k[\Delta_m])/ \Theta \Omega(k[\Delta_m]) \cong I/\Theta I \cong I/(I \cap \Theta).$$  

Consider the commutative diagram,
  
  $$\begin{array}{ccccccccc}
   0 & \to & I/(I \cap \Theta)_i &\to& k(\Delta)_i &\to& k(\Sigma)_i &\to &0 \\
    & & & & & & &\\
   & & \downarrow \omega^{d-2i}& &\ \ \downarrow \omega^{d-2i}& &\ \ \downarrow \omega^{d-2i}& & \\
     & & & & & & &\\
      0 & \to & I/(I \cap \Theta)_{d-i} &\to& k(\Delta)_{d-i}&\to& k(\Sigma)_{d-i} &\to &0
   \end{array}$$

  \noindent  As $\Omega(k[\Delta_m])/\Theta \Omega(k[\Delta_m]) \cong \Hom_k (k(\Delta_m),k), \ 
I/(I \cap \Theta)$ must also be isomorphic to $\Hom_k (k(\Delta_m),k).$  By Theorem \ref{surj balls}, multiplication  $\omega^{d-2i}:k(\Delta_m)_i \to k(\Delta_m)_{d-i}$ is a surjection for the nonempty Zariski  open set $G(\partial P_m),$ where $P_m$ is a simplicial $d$-polytope such that $\Delta_m \subset \partial P_m.$  As the left-hand vertical arrow is the $k$-dual for this map, it must be an injection for pairs $(\omega, \Theta) \in G(\partial P_m).$   The induction hypothesis provides another nonempty Zariski open subset of pairs  such that the right-hand  vertical arrow is an injection.  The intersection of these two sets is a nonempty Zariski open subset such that the middle vertical arrow is an injection. 
  \end{proof} 
  
We note that if at some point in the  future Theorem \ref{g thm} is extended to a more general class of homology spheres, say $\SS,$ then it would be possible to define $\SS$ ear decompositions.  In that case, the above proof would still be valid.

\begin{cor} \label{main enumerate}
  Let $\Delta$ be a $(d-1)$-dimensional simplicial complex with a convex ear decomposition.  
  \begin{enumerate}
    \item[a.] \label{Chari}
       If $i \le \lfloor d/2 \rfloor,$ then $h_i \le h_{d-i}$ and $h_i \le h_{i+1}.$
    \item[b.]
       If $g_i = h_i - h_{i-1},$ then $(g_0, g_1, \dots, g_{\lceil d/2 \rceil})$ is an M-vector.
   \end{enumerate}
\end{cor}

\begin{proof}
Let $\omega$ be a $g$-element for $k(\Delta)$  and let $i \le \lfloor d/2 \rfloor.$  As multiplication $\omega^{d-2i}:k(\Delta)_i \to k(\Delta)_{d-i}$ is an injection, $h_i \le h_{d-i}$ and multiplication $\omega:k(\Delta)_i \to k(\Delta)_{i+1}$ is also an injection.  Hence $h_i \le h_{i+1}.$  To see that $(g_0, g_1, \dots, g_{\lceil d/2 \rceil})$ is an M-vector we simply note that for $i \le \lceil d/2 \rceil,$ $g_i = \dim_k (k(\Delta)/<\omega>)_i.$
\end{proof}

\noindent The inequalities in (a) are originally due to Chari \cite{Ch}.

Theorem \ref{main} was first proved for independence complexes of rationally represented matroids by Hausel and Sturmfels \cite{HS}.  They used the theory of hyperk\"{a}hler toric varities.   Swartz proved this for all matroids \cite{Sw3}.  In \cite{Ha} Hausel presents a proof for all matroids which is based on ideas from both papers.  The above proof is an extension of this idea to spaces with a convex ear decomposition.

\end{section}

\begin{section}{Higher Cohen-Macaulay connectivity} \label{kcm}

In the previous section we said that spaces with a convex ear decomposition are doubly Cohen-Macaulay.  A Cohen-Macaulay complex $\Delta$ is {\it doubly Cohen-Macaulay} if for all vertices $v$ in $\Delta,$ the dimension of $\Delta - v$ equals the dimension of $\Delta,$ and $\Delta-v$ is still Cohen-Macaulay.  More generally, $\Delta$ is {\it $q$-CM} if $\Delta$ is Cohen-Macaulay and for every subset $A$ of vertices of $\Delta$ with $|A| < q, \ \dim (\Delta-A )= \dim \Delta$ and $\Delta-A$ is still Cohen-Macaulay.  The maximum $q$ such that $\Delta$ is $q$-CM is the {\it CM-connectivity} of $\Delta.$  

\begin{thm}
 If $\Delta$ has a convex ear decomposition, then $\Delta$ is doubly Cohen-Macaulay.
\end{thm}

As is conjecturally  the case with Theorem \ref{main}, this result extends to a more general construction using homology spheres and balls since we only use the homological properties of balls and spheres in the proof.  Here, a {\it homology sphere} has the homology of a sphere and every link has the homology of a sphere of the appropriate dimension. Homology spheres are also called {\it Gorenstein*} complexes.  A {\it homology ball} is homologically acyclic and every link is either a homology ball or a homology sphere of the appropriate dimension.  Furthermore, the faces whose links are homology balls form a subcomplex, the boundary,  which is a homology sphere of one lower dimension.  Removing a vertex from a homology sphere leaves a homology ball with boundary the link of the vertex in the original homology sphere.  

\begin{proof}
    Let $v$ be a vertex in $\Delta$ and let  $\sigma$ be a face of $\Delta -v.$  We must show  $\lk_{\Delta-v} \sigma$ is Cohen-Macaulay.  The link of $\sigma$ in $\Delta -v$ is the union of the links of  $\sigma$ in $\tilde{\Delta}_1-v, \dots \tilde{\Delta}_t-v,$ where the $\tilde{\Delta}_j$ is a sequential renumbering of all the $\Delta_i$ which contain $\sigma.$  The proof is by induction on $t.$  If $t=1,$ then the link is either a homology ball or homology sphere, depending on whether or not $v$ is in the link of $\sigma$ in $\tilde{\Delta}_1.$  In either case the link is Cohen-Macaulay.

When $j > 1, \ \sigma$ is on the boundary of $\tilde{\Delta}_j,$  hence its link in $\tilde{\Delta}_j$ is a homology ball.  Let $Y=\cup^{t-1}_{j=1} \lk_{\tilde{\Delta}_j -v} \sigma$ and $Z=\lk _{\tilde{\Delta}_t -v} \sigma.$ For the induction step there are three possibilities.  In each case, Mayer-Vietoris arguments are sufficient to compute the homology of $Y \cap Z, Z$ and see that their union, $\lk_{\Delta-v} \sigma,$ is Cohen-Macaulay. 

\begin{enumerate}
  \item  $v \notin \lk_{\tilde{\Delta}_t} \sigma.$  Then $Z$ is a $(d-1-|\sigma|)$-homology ball and $Y \cap Z$ is a $d-2-|\sigma|$-homology sphere.
  \item  $v$ is an interior point of $\lk_{\tilde{\Delta}_t} \sigma.$  Then $Y \cap Z$ and $Z$ have the homology of a $(d-2-|\sigma|)$-homology sphere. For $Z,$ this follows from the Mayer-Vietoris sequence for the homology ball $\lk_{\tilde{\Delta}_t} \sigma$ written as the union of $Z$ and the closed star of $v$ in $\lk_{\tilde{\Delta}_t} \sigma$. In addition, the inclusion map $Y \cap Z \to Z$ is an isomorphism in homology.
  \item  $v$ is a boundary point of $\lk_{\tilde{\Delta}_t} \sigma.$  Now $Y \cap Z$ and $Z$  are homologically acyclic.
  \end{enumerate}
 \end{proof}

Doubly Cohen-Macualay complexes and  $q$-CM complexes were introduced by Baclawski \cite{Ba2}. Spheres, and more generally homology spheres, are $2$-CM, but balls are not.   Baclawski proved that the order complex of a  semimodular poset $P$ is $2$-CM if and only if $P$ is a geometric lattice. Furthermore, if $P$ is a geometric lattice, then its order complex  is $q$-CM if and only if every line of $P$ has at least $q$ atoms.  In his study of buildings and Coxeter complexes Bj\"{o}rner proved that finite buildings are $2$-CM and any finite building associated to $GF(q)$ is $(q+1)$-CM \cite{Bj3}.  Welker showed that order complexes of supersolvable lattices are $2$-CM if and only if the M\"{o}bius function is nonzero on every interval \cite{We}.  Since all of the above examples of $2$-CM complexes either have a convex ear decomposition, and hence satisfy the conclusion of Theorem \ref{main}, or are conjectured to satisfy Theorem \ref{main}, it seems natural to ask the following question which was also suggested by  Bj\"orner.

\begin{prob} \label{prob}
  Do all $2$-CM complexes satisfy the conclusion of Theorem \ref{main}?
\end{prob}

 Nevo has shown that for  $2$-CM complexes and $d \ge 3, \omega: k(\Delta)_1 \to k(\Delta)_2$ is  injective for a generic set of pairs $(\omega, \Theta)$  \cite{Ne}. An affirmative answer to this question would also show that all homology spheres satisfy the $g$-theorem \cite[Conjecture II.6.2]{St}.  The face ring of any 2-CM complex is a level ring (see \cite[pg. 94]{St} for a discussion).  Combined with \cite{Hi}, this implies that for any 2-CM complex
     $$h_0 + h_1 + \dots + h_i \le h_d + h_{d-1} + \dots + h_{d-i}.$$
 
Another property shared by all $(d-1)$-dimensional 2-CM complexes which comes from the fact that their face rings are level, is that the reversed $h$-vector, $(h_d, \dots, h_0)$ is a sum of $h_d \ M$-vectors \cite{St1}.    For spaces with a convex ear decomposition, this fact expresses itself in the following formula \cite{Ch}.

\begin{equation} \label{comph}
h_{d-i}(\Delta) = h_i(\Delta_1) + \dots + h_i(\Delta_m) .
\end{equation}

 When $\Delta$ has a convex ear decomposition, $h_{d-i} - h_i$ is nonnegative for $i \le d/2.$ So it seems natural to consider the following complementary $h$-vector of $\Delta.$
 
 \begin{defn}
   The {\bf complementary $h$-vector} of $\Delta$ is
     $$  \bar{h} = (h_d - h_0, h_{d-1}- h_1, \dots, h_{d- \lfloor d/2 \rfloor} - h_{d- \lceil d/2 \rceil}).$$
   \end{defn}
   
   For any homology ball $\Delta$, $h_i(\Delta) - h_{d-i}(\Delta) = g_i(\partial \Delta)$ (see, for instance, \cite{MW}).  If $(\Delta_1, \dots \Delta_m)$ is a convex ear decomposition for $\Delta, $ then equation (\ref{comph}) implies
   
     $$h_{d-i}(\Delta) - h_i(\Delta) = \displaystyle\sum^m_{j=2}h_i ( \Delta_j)- h_{d-i}(\Delta_j) =  \displaystyle\sum^m_{j=2} g_i(\partial \Delta_j).$$  
     
     \noindent   Thus, if the boundaries of the ears are known to be combinatorially equivalent to the boundaries of $(d-1)$-polytopes, then by Theorem \ref{g thm} $\bar{h}$ is the sum of $h_d - 1$ M-vectors.
     
     \begin{prop}  If $\Delta$ has a PS-ear decomposition, then $\bar{h}$ is a sum of $h_d - 1$ M-vectors.  
     \end{prop}
     
     \begin{proof}
       The boundary of each ear in a PS-ear decomposition is the join of the boundaries of simplices.
      \end{proof}
      
      In order to analyze $\bar{h}$ when $\Delta$ is a finite building we use von Heydebreck's  convex ear decomposition.
      
      \begin{prop}\cite[Lemma 3.4]{He}
        Let $\Delta$ be a finite building of type $(W,S)$.  Then $\Delta$ has a convex ear decomposition such that each ear is isomorphic to $\Sigma(W,S) \cap H_1^+\cap \dots \cap H_t^+,$ where the $H_i^+$'s are closed  half-spaces of distinct reflecting hyperplanes  of $W.$
       \end{prop}
       
    \begin{lem}
         Let $\A = \{H_1, \dots, H_s\}$ be an essential arrangement of hyperplanes in $\R^d.$  Let $P$ be any $d$-polytope whose face fan is the fan of $\A.$ Let $H^+_{i_1}, \dots, H^+_{i_t}$ be  closed half-spaces of distinct hyperplanes in $\A.$ If $B = \partial P \cap H^+_{i_1} \cap \dots \cap H^+_{i_t}$ is nonempty, then $\partial B$ is combinatorially equivalent to the boundary of a $(d-1)$-polytope.      \end{lem}  
       
   \begin{proof}
   For notational convenience we renumber the hyperplanes so that $H_{i_j}=H_j.$  So, $\A=\{H_1,\dots,H_t,H_{t+1},\dots,H_s\}$ and $B = \partial P \cap H^+_1 \cap \dots \cap H^+_t.$
   
        By a familiar vertex figure argument, it is sufficient to find $P^\prime,$ a $d$-polytope whose face fan is the same as the face fan of $P,$ and a point $y \in \R^d$ such that the facets of $P^\prime$ that can be seen from $y$ are precisely those in $B^\prime = \partial P^\prime \cap H^+_1 \cap \dots \cap H^+_t$.
        
         Let $z_1, \dots , z_s$ be nonzero vectors such that $z_i$ is orthogonal to $H_i,$  and for $i \le t, \ z_i \in H^+_i.$  Let $Z$ be the zonotope $[-z_1, z_1] + \dots +[-z_s, z_s],$ where $+$ is Minkowski sum. For an $s$-tuple $(\varepsilon_1, \dots, \varepsilon_s), \ \varepsilon_i = \pm 1,$ let   $x_\varepsilon = \displaystyle\sum^s_{i=1} \varepsilon_i z_i.$   While some of the $x_\varepsilon$ may be interior points of $Z,$ all of the vertices of $Z$ are equal to some $x_\varepsilon.$ 
         
          Set $P^\prime$ to be $Z^\star,$ the polar of $Z.$  The face fan of $P^\prime$ is the fan of $\A$ and the facets of $P^\prime$ are of the form $F_\varepsilon = P^\prime \cap H_\varepsilon,$ where $H_\varepsilon = \{w \in \R^d: x_\varepsilon \cdot w = 1, x_\varepsilon \mbox{ a vertex of } Z \}.$ If $x_\varepsilon$ is a vertex of $Z,$ then $x_\varepsilon \cdot w \le 1$ for all $w \in P^\prime.$  Therefore, a facet $F_\varepsilon$ of $P^\prime$ is visible from $y \in \R^d$ if and only if $x_\varepsilon \cdot y > 1.$  
          
           The facets of $B^\prime = P^\prime \cap H^+_1 \cap \dots \cap H^+_t$ are  those facets $F_\varepsilon$  with $\varepsilon_i = +1$ for $i \le t.$  Since $B$ is nonempty we can choose $y \in R^d$ so that $z_i \cdot y > 0$ for $i \le t.$  By rescaling all the $z_i$ by positive scalars, we can choose $z_1, \dots, z_s$ and $y$ so that $ z_i \cdot y = 1$ for $i \le t,$ and $\displaystyle\sum^s_{i=t+1} |z_i \cdot y| < 1/100.$  Now it is easy to see that for some $\delta > 0$ the only facets of $P^\prime$ visible from $(1+ \delta) y $ will be the facets of $B^\prime.$  
         \end{proof}
         
   \begin{cor}
     If $\Delta$ is a finite building, then $\bar{h}$ is the sum of $h_d -1$ M-vectors.
   \end{cor}
   
   If the g-theorem was known to hold for the boundary of every ball which is a full dimensional subcomplex of the boundary of a simplicial polytope, then $\bar{h}$ of any $d$-dimensional space with a convex ear decomposition would be the sum of $h_d-1$ M-vectors.
   
   \begin{prob}
     If $\Delta$ is 2-CM, is $\bar{h}(\Delta)$ the sum of $h_d -1$ M-vectors?
   \end{prob}

Bj\"{o}rner has observed that for $q >> 0,$ equation (\ref{GF(q) flagh})  implies that finite buildings associated to $GF(q)$ must satisfy $h_B < h_A$ for $B \subset A$ \cite{Bj4}.  This means that for these spaces high CM connectivity implies that the $h$-vector is increasing.  This is part of a general phenomenon involving complexes with large links.  Since the removal of $q-1$ vertices from a $q$-CM complex $\Delta$ leaves a pure complex, the link of every nonfacet  of $\Delta$ must contain at least $q$ vertices.  
A pure complex with large links must have an increasing $h$-vector.   To prove this we use an extension of $\tilde{h},$  the short simplicial $h$-vector introduced in \cite{HN} as a simplicial analogue of the short cubical $h$-vector in \cite{Ad}.

\begin{defn}
  $ \tilde{h}^{(m)}_i(\Delta) = \displaystyle\sum_{\stackrel{\sigma \in \Delta}{|\sigma|=m}} h_i(\lk \sigma).$
\end{defn}

\noindent For example, $\Th^{(0)}$ is the usual $h$-vector and  $\Th^{(1)}$ is the short simplicial $h$-vector defined in \cite{HN}.

\begin{prop}
  Let $\Delta$ be a pure $(d-1)$-dimensional complex. For $m \le d-1$ and $i \le d-m,$ 
  \begin{equation} \label{ss h}
  (m+1)\  \Th^{(m+1)}_{i-1} = i \Th^{(m)}_i + (d-m-i+1) \Th^{(m)}_{i-1}.
  \end{equation}
\end{prop}

\begin{proof}
  The case $m=0$ is  \cite[Proposition 2.3]{Sw5}.  For larger $m,$
  $$
  \Th^{(m)}_i   =  \displaystyle\sum_{|\sigma|=m} h_i(\lk \sigma)$$
  $$\begin{array}{cl}
  = & \frac{1}{i} \left[ \displaystyle\sum_{|\sigma|=m} \{\Th^{(1)}_{i-1} (\lk \sigma) - (d-m-i+1) h_{i-1}(\lk \sigma) \} \right]\\
   = & \frac{1}{i} \left[ \displaystyle\sum_{|\sigma|=m} \left( \left\{ \displaystyle\sum_{v \in \lk \sigma} h_{i-1}(\lk (\sigma \cup \{v\})) \right\} - (d-m-i+1) h_{i-1}(\lk \sigma) \right) \right] \\
   = & \frac{1}{i} \left[ (m+1) \ h^{(m+1)}_{i-1} - (d-m-i+1) h^{(m)}_{i-1} \right].
   \end{array}
   $$
\end{proof}

\begin{thm}
  Fix $d$ and $i \le d.$  There exists $q(i,d)$ such that if $\Delta$ is a  pure $(d-1)$-dimensional complex and the link of every $(i-2)$-dimensional face of $\Delta$ has at least $q(i,d)$ vertices, then 
  
\begin{equation} \label{increasing h}
h_0 < \dots < h_i.
\end{equation}
\end{thm}

\begin{proof}
 The proof is by induction on $i$ and $d$, with $i=1$ being trivial.  We can assume that $q(i,d) \le q(i^\prime, d^\prime)$ whenever $i \le i^\prime$ and $d \le d^\prime.$  Suppose that $i \ge 2$ and $q(i-1,d)$ satisfies the theorem.  For the induction step we only need to find $q$ so that $h_{i-1} < h_i$ whenever $\Delta$ is pure, $(d-1)$-dimensional and the link of every $(i-2)$-dimensional face has at least $q$ vertices.  Indeed, given such a $q,\   q(i,d) = \max \{q, q(i-1,d)\}$ satisfies the theorem.
 
Let $q$ be the minimum number of vertices in the link of an $(i-2)$-dimensional face of $\Delta.$  In order to show that for $q$ sufficiently large $h_{i-1} < h_i$ we argue by contradiction.  So, suppose $h_{i-1} \ge h_i.$ First we estimate $\Th^{(j)}_{i-j}$ using (\ref{ss h}) and $h_{i-1} \ge h_i.$ For $0 \le j \le i-1,$

\begin{equation} \label{above}
\Th^{(j)}_{i-j} \le p_{i,j}(d) h_{i-1},
\end{equation}

\noindent where $p_{i,j}(d)$ is a degree $j$ polymomial.  When $j=0, \ p_{i,0} = 1.$  For higher $j,$ (\ref{above}) is proved by induction using $\Th^{(j)}_{i-j-1} < \Th^{(j)}_{i-j}$ and

$$\Th^{(j+1)}_{i-j-1} = \frac{1}{(j+1)}\ \left[ (i-j) \Th^{(j)}_{i-j} + (d-i+1) \Th^{(j)}_{i-j-1} \right].$$

We can also estimate $\Th^{(i-1)}_1$ using the induction hypothesis and $q$.  

$$ \Th^{(i-1)}_1 = \displaystyle\sum_{|\sigma|=i-1} h_1(\lk \sigma) \ge (q-(d-i+1)) f_{i-1} > (q-(d-i+1)) h_{i-1}.$$

\noindent  The last inequality follows from (\ref{f by h}) and the fact that $h_0, \dots, h_{i-1}$ are positive for sufficiently large $q$ by the induction hypothesis on $i.$ Putting these estimates together,

$$(q-(d-i+1)) h_{i-1} < p_{i,i-1}(d) h_{i-1}.$$

\noindent  Hence, $q$ is bounded.  
\end{proof}

The above proof gives $q(i,d)$ as  a polynomial of degree $i-1$ in $d.$  Can this be improved?

\begin{prob}  What are the minimum values of $q(i,d)$? \end{prob}

The dependence on $d$ in the above theorem is essential.  Indeed, $h_2$ of two   simplices connected at one vertex is always negative once $d > 2.$  What if we impose the additional condition that $\Delta$ is Cohen-Macaulay? An affirmative answer to Problem \ref{prob} would imply that links of size $q(i,2i)$ would be sufficient to imply at least inequality (as opposed to strict inequality) in (\ref{increasing h})  for $2$-CM complexes.   For  balanced $q$-CM complexes it is possible to remove the dependence on $d.$  

\begin{thm}  \label{colored}
  Let $(\Delta, \psi)$ be a balanced $q(i,i)$-CM complex. Then for any $B \subset A$ with $|A| \le i,$
  $$ h_B < h_A.$$  
\end{thm}

\begin{proof}
 Our first observation is that for any $A \subseteq S, \ \Delta_A$ is also $q(i,i)$-CM.  Removing $q(i,i)-1$ vertices from $\Delta_A$ is the same as removing $q(i,i)-1$ vertices from $\Delta$ and then restricting to $\Delta_A.$ The combination of these operations preserve dimension and the CM property \cite[Theorem III.4.5]{St}.  Hence,

$$h_B(\Delta) = h_B(\Delta_A) \le h_{|B|}(\Delta_A) < h_{|A|}(\Delta_A) = h_A(\Delta).$$

\end{proof}

\begin{cor}
  If $\Delta$ is a $q(i,i)$-CM balanced complex and $d \ge i +1,$ then
  $$ h_0 < \dots < h_i.$$
\end{cor}

\begin{proof}

In order to show that $h_{j-1} < h_j$ for every $j \le i,$ fix $j \le i.$ Let $\phi$ be any map from $(j-1)$-subsets of $S$ to $j$-subsets of $S$ such that $B \subseteq \phi(B)$ for all $B.$ 

$$h_{j-1}(\Delta) = \sum_{|B| = j-1} h_B(\Delta) = \sum_{|B|=j-1} h_B(\Delta_{\phi(B)}) =\sum_{|A|=j}\  \sum_{B \in \phi^{-1}(A)} h_B (\Delta_A)$$
$$ \le \sum_{|A|=j} h_{j-1}(\Delta_A)  < \sum_{|A|=j} h_j(\Delta_A) = \sum_{|A|=j} h_A(\Delta) = h_j(\Delta).$$
\end{proof}

\end{section}

{\it Acknowledgements:}  Kai-Uwe Bux pointed out the relevance of \cite{He}. Anders Bj\"orner suggested the question of when $h_B < h_A$ and the importance of equation (\ref{GF(q) flagh}).

 \bibliographystyle{plain}

\end{document}